\documentclass[letterpaper, 10 pt, conference]{ieeeconf}  

\IEEEoverridecommandlockouts                              

\usepackage{mathptmx} 
\usepackage{times} 
\usepackage{amssymb}  
\usepackage{graphicx}
\usepackage{amsmath, nccmath}
\usepackage{tabularx,ragged2e,booktabs,caption}
\usepackage{geometry}
\usepackage{gensymb}

\title{\LARGE \bf
Indirect Adaptive Fuzzy Model Predictive Control of a Rotational Inverted Pendulum
}

\author{Roja Eini$^{1}$ and Sherif Abdelwahed$^{2}$
	\thanks{$^{1}$Roja Eini is PhD student of Department of Electrical Engineering,
		Virginia Commonwealth University, Richmond, VA, USA
		{\tt\small einir@vcu.edu}}%
	\thanks{$^{2}$Sherif Abdelwahed is faculty of Department of Electrical Engineering,
		Virginia Commonwealth University, Richmond, VA, USA
		{\tt\small sabdelwahed@vcu.edu}}%
}

\begin{document}

\maketitle
\thispagestyle{empty}
\pagestyle{empty}

\begin{abstract}
This paper introduces an indirect adaptive fuzzy model predictive control strategy for a nonlinear rotational inverted pendulum with model uncertainties. In the first stage, a nonlinear prediction model is provided based on the fuzzy sets, and the model parameters are tuned through the adaption rules. In the second stage, the model predictive controller is designed based on the predicted inputs and outputs of the system. The control objective is to track the desired outputs with minimum error and to maintain closed-loop stability based on the Lyapunov theorem. Combining the adaptive Mamdani fuzzy model with the model predictive control method is proposed for the first time for the nonlinear inverted pendulum. Moreover, the proposed approach considers the disturbances predictions as part of the system inputs which have not been considered in the previous related works. Thus, more accurate predictions resistant to the parameters variations enhance the system performance using the proposed approach. A classical model predictive controller is also applied to the plant, and the results of the proposed strategy are compared with the results from the classical approach. Results proved that the proposed algorithm improves the control performance significantly with guaranteed stability and excellent tracking. 
\end{abstract}  


\section{INTRODUCTION}

Model Predictive Control (MPC) strategy has been recently known as a powerful control method for industrial applications, especially for highly nonlinear systems with uncertainties and constraints [1]. The most prominent part of a model-based predictive approach is to find a prediction model to approximate the system's future input, output, and state signals. Defining an accurate prediction model is a difficult problem in predictive control algorithms since in practical applications the certain mathematical model of the system is not accessible. Moreover, when the system is nonlinear, the optimization problem becomes non-convex, and there is no benefit from the use of standard prediction models in the predictive algorithms. Researchers have used various system modeling approaches such as neural networks and fuzzy systems to obtain accurate prediction model for predictive algorithms [2-7]; however, most of the related works in this area are applied to either a linearized plant or a nonlinear plant without considering its uncertainties.    \\ \indent
Fuzzy systems are designed based on previous experiences or specific knowledge of the system. Indeed, if-then fuzzy rules are very applicable in system modeling, especially in modeling unknown nonlinear systems with uncertainties. Among fuzzy systems, adaptive fuzzy systems are more powerful because their model parameters are robust to uncertainties and disturbances. Adaptive fuzzy systems are categorized into two frameworks: indirect adaptive fuzzy systems, and direct adaptive fuzzy systems. In designing an indirect adaptive fuzzy system, a fuzzy system is first constructed, and then its parameters are regulated based on the adaption rules.  \\ \indent
This paper presents how a controller based on model predictive control theory can be developed based on an indirect adaptive fuzzy model. Nonlinear rotational inverted pendulum system is highly nonlinear with time-varying parameters and model uncertainties. A nonlinear Mamdani fuzzy model is considered as the prediction model. The fuzzy system's parameters are tuned based on the adaption rules such that the Lyapunov stability criterion is maintained and the system tracking error is minimized. Combination of the fuzzy modeling approach and model predictive control method is applied for the first time to the nonlinear rotational inverted pendulum benchmark problem. To evaluate the effectiveness of the proposed control methodology, a classical model predictive controller is also designed for the system, and the results of the two controllers are compared.   \\ \indent
The rest of the paper is organized as follows. The system model and its parameters are described in section II. Section III introduces the classical MPC approach and the proposed indirect adaptive fuzzy MPC methodology. In section IV, the simulation results are shown, and the last section provides the conclusions.
 
\section{SYSTEM DEFINITION}

\begin{figure}[thpb]
	\begin{center}
		\framebox{
			\parbox{2in}
			{
				\centering
				\includegraphics[height=4cm, width=5cm]{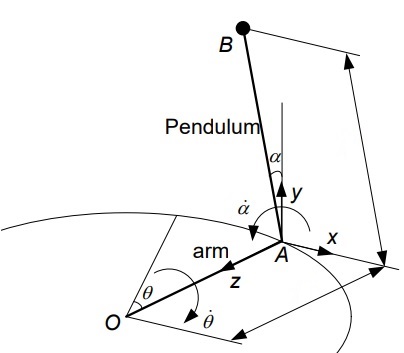}
			}
		}
		\caption{Simple model of a rotational inverted pendulum}
		\label{fig:pendulum}
	\end{center}
\end{figure}

The nonlinear rotational inverted pendulum system is known as an interesting and applicable system in control of missile launchers, pendubots, segways, and earthquake resistants. The objective of controlling an inverted pendulum is to balance pendulum in its upright position. Fig. \ref{fig:pendulum} shows a simple model of a rotational inverted pendulum consists of a vertical arm, a horizontal rotating arm, a gear chain, and a DC motor. A rotary encoder is also attached to the shaft to feedback the angle of the arms and to measure the pendulum's motion [8].  \\ \indent
The mathematical model of a rotational inverted pendulum is as (\ref{eq:pendulum}). The model consists of four states $\theta$, $\theta^\cdot$, $\alpha$, and $\alpha^\cdot$. The two states $\alpha$ and $\theta$ are the pendulum's vertical angle and the pendulum's horizontal angle respectively, which can be measured. $\theta^\cdot$ and $\alpha^\cdot$ can also be attained by the first-order derivatives of $\theta$ and $\alpha$ respectively. Moreover, $u$ is the pendulum's input signal. Besides, $m_1$, $a_p$, $k_p$, $k_1$, $g$, $l_1$, $J_1$, and $c_1$ are the parameters of the plant as Table \ref{tab:parameters}. \\ \indent
\begin{fleqn}
	\begin{equation}  
	\begin{aligned}[b]
	& \ddot \theta=a_p \dot \theta + k_pu\\
	& \ddot \alpha=-\frac{c_1}{J_1} \dot \alpha+\frac{m_1gl_1}{J_1}sin\alpha + \frac{k_1}{J_1} \ddot \theta
	\end{aligned}
	\label{eq:pendulum}
	\end{equation}
\end{fleqn}
\\ \indent
\begin {table}
\caption {system parameters}
\begin{center}
	\noindent
	\newlength{\myheight}
	\setlength{\myheight}{0.5cm}
	\begin{tabular}{|l|c|}  
		\hline  \hline
		\parbox[c][\myheight][c]{0cm}{}   $m_1$   &  0.0861 kg \\  
		\parbox[c][\myheight][c]{0cm}{}    $k_1$  &  0.0019
		\\   
		\parbox[c][\myheight][c]{0cm}{}    $a_p$  & 33.04 
		\\  
		\parbox[c][\myheight][c]{0cm}{}    $J_1$  & 0.0010
		\\ 
		\parbox[c][\myheight][c]{0cm}{}    $g$    &  9.8066   
		\\
		\parbox[c][\myheight][c]{0cm}{}    $l_1$  &  0.113 m  
		\\
		\parbox[c][\myheight][c]{0cm}{}    $c_1$  &  0.0029  
		\\
		\parbox[c][\myheight][c]{0cm}{}    $k_p$  &  74.89  \\
		\hline \hline 
	\end{tabular}
	\label{tab:parameters} 
\end{center}
\end {table} 
The sate space representation of the pendulum system is represented as (\ref{eq:statespace}); where $x_1=\theta$, $x_2=\dot \theta$, $x_3=\alpha$, and $x_4=\dot \alpha$. \\ 

\begin{fleqn}
	\begin{equation}  
	\begin{aligned}[b]
	&\dot x_1 = x_2\\
	&\dot x_2 = a_1x_2 + b_1u\\
	&\dot x_3 = x_4\\
	&\dot x_4 = a_2x_2 + a_3 sinx_3 + a_4 x_4 + b_2u
	\end{aligned}
	\label{eq:statespace}
	\end{equation}
\end{fleqn}  \\ \indent
Additionally, $a_1=-a_p$, $a_2=-\frac{k_1a_p}{J_1}$, $a_3=\frac{m_1gl_1}{J_1}$, $a_4=-\frac{c_1}{J_1}$, $b_1=k_p$, and $b_2=\frac{k_1k_p}{J_1}$. Hence, the aim is to balance the pendulum's output $y=x_3$.  \\

\section{THE MODEL PREDICTIVE CONTROL APPROACH}

\begin{figure}[thpb]
	\begin{center}
		\framebox{\parbox{2.5in}{
				\centering
				\includegraphics[height=4cm, width=5cm]{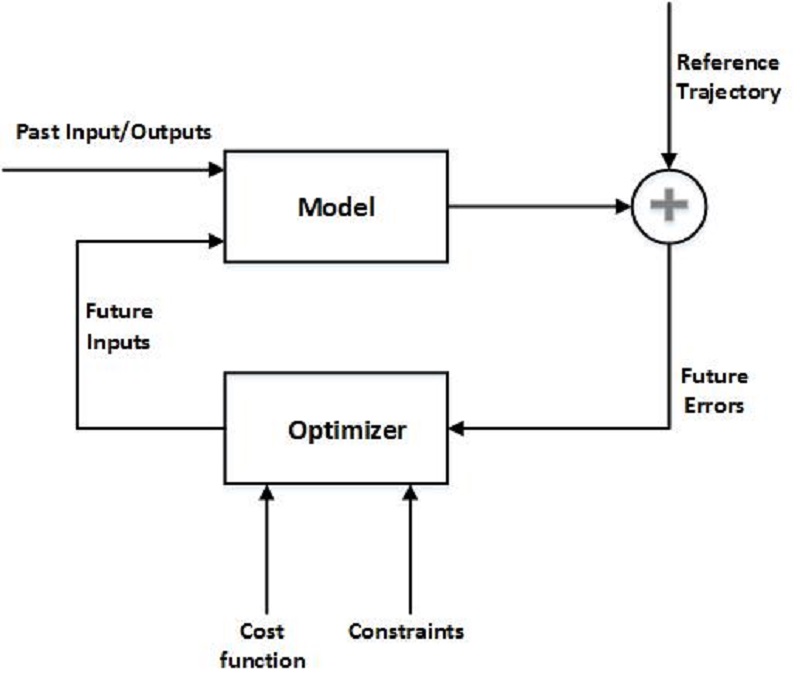}
		}}
		\caption{Simple structure of a model predictive control strategy}
		\label{MPC}
	\end{center}
\end{figure}

Since rotational inverted pendulum is a highly nonlinear system subject to the constraints and model uncertainties, applying model-based predictive approaches to control it is considered as a reliable strategy.  Model predictive control allows for the real-time control of a nonlinear inverted pendulum through predicting the system states, outputs, disturbances, and inputs. Note that the predictive models are necessary to estimate the signals throughout the prediction horizon [1]. Fig. \ref{MPC} shows a simple structure of model predictive control.
\\ \indent
Model predictive control methods use different prediction models. Herein, the classical model predictive control method is described in the first subsection, and the proposed model predictive control approach based on indirect adaptive fuzzy system is introduced in the next one. \\ 

\subsection{THE CLASSICAL MODEL PREDICTIVE CONTROL}
The nonlinear system of (\ref{eq:statespace}) can be stated as a simple nonlinear model of (\ref{eq:general_nonlinear}). \\
\begin{fleqn}
	\begin{equation}  
	\begin{aligned}[b]
	&x(k+1)=f_k(x(k),u(k))
	\end{aligned}
	\label{eq:general_nonlinear}
	\end{equation}
\end{fleqn}
where $x(k)$ and $u(k)$ denote the state and the input vectors at time $k$ respectively. $f_k$ is the mapping function from the current states and inputs to the next state and input values. Therefore, the system model from the first step to the last step in the horizon is defined as (\ref{eq:general_nonlinear_descriptive}). Note that the control horizon and the states horizon are defined as $N_c$ and $N_p$ respectively. Thus, the input sequence $U$ would be (\ref{eq:u}).  \\
\begin{fleqn}
	\begin{equation}  
	\begin{aligned}[b]
	&x(2)=f_1(x(1),u(1))  \\
	&  \vdots  \\
	&x(N_p+1)=f_1(x(N_p),u(N_c))
	\end{aligned}
	\label{eq:general_nonlinear_descriptive}
	\end{equation}
\end{fleqn} 
\begin{fleqn}
	\begin{equation}  
	\begin{aligned}[b]
	&U=(u(1),\,\, \hdots,\,\, u(N_c)) 
	\end{aligned}
	\label{eq:u}
	\end{equation}
\end{fleqn}  \\ \indent 
The cost function for the classical predictive control algorithm is defined as (\ref{cost}). The  control objective is to minimize the tracking error toward the desired trajectory with minimum control effort. Thus, the cost function includes the quadratic form of the control signal $u$ and the state vector $x$. \\
\begin{fleqn}
	\begin{equation}  
	\begin{aligned}[b]
	&J(k)=\sum\limits_{p=1}^{K_p} {\lVert{\hat{x}(k_j+p|k)\lVert}_Q }^2 
		&+ \sum\limits_{p=1}^{K_c} {\lVert{\hat{u}(k_u+p|k)\lVert}_R }^2 
	\end{aligned}
	\label{cost}
	\end{equation}
\end{fleqn}  \\ 
where $Q$ and $R$ are the positive definite weighting matrices of the states and the control signals respectively. Furthermore, $\hat x(k)$ and $\hat u(k)$ are the predicted state and predicted input signals at time step $k$ respectively. \\ \indent 
The model constraints are the predicted state and input equations attained from the prediction model (\ref{prediction}). In the classical predictive control method, the prediction model is a simple impulse response model which also includes the disturbance signal [1]. \\
\begin{fleqn}
	\begin{equation}  
	\begin{aligned}[b]
	&\hat x(k+1)=f(\hat x(k),\hat u(k),d(k))
	\end{aligned}
	\label{prediction}
	\end{equation}
\end{fleqn}  \\ 
where $d(k)$ is the system disturbance at time $k$. \\ \indent
Therefore, the cost function $J(k)$ is minimized in each step of the model predictive control algorithm (over the whole control and state horizon). In each step only the first control law is applied to the system and the procedure repeats in the next step. \\ \indent
Thus, the optimized state is defined as (\ref{state_optimum}). \\
\begin{fleqn}
	\begin{equation}  
	\begin{aligned}[b]
	&(x_{opt}, u_{opt})=arg\{\,min \,\, J(k)\}, \,\,\,\,c(x,u)\leq 0
	\end{aligned}
	\label{state_optimum}
	\end{equation}
\end{fleqn}  \\ 
where $c$ is the inequality function, and $x_{opt}$ and $u_{opt}$ are the optimum variables from solving the optimization problem. The Lagrangian function $L(x,\lambda)$ is defined in (\ref{Lagrange}). \\
\begin{fleqn}
	\begin{equation}  
	\begin{aligned}[b]
	&L(x,\lambda)=J(k)+\lambda^Tc(x,u)
	\end{aligned}
	\label{Lagrange}
	\end{equation}
\end{fleqn}  \\ \indent
Using the fmincon optimization tool in MATLAB, the optimum values for the search direction variable $p$ and Lagrange multiplier $\lambda$ is attained as (\ref{pandlambdaoptimum}) [6]. 
\begin{fleqn}
	\begin{equation}  
	\begin{aligned}[b]
	&(p_{opt},\lambda_{opt})=arg\{\,min \,\, p^TL_x(x,\lambda)+\frac{1}{2}P^TL_{xx}(x,\lambda)P\}
	\end{aligned}
	\label{pandlambdaoptimum}
	\end{equation}
\end{fleqn}  \\ 
where $p_{opt}$ and $\lambda_{opt}$ are the optimum values of $p$ and $\lambda$ respectively. \\ \indent
Hence, the classical model predictive control algorithm is as follows [9]: \\
\begin{itemize}
	\item at time $k=0$, determine the state value $x(0)$ and get the input value $u(0)$ by solving the optimization problem. 
	\item at time $k>0$, get the state, input, and disturbance predictions from (\ref{prediction}) (with $K_c$ as the disturbance and input horizon, and $K_p$ as the state horizon).
	\item at time $k>0$, solve the optimization problem (\ref{cost}) to get the optimum input signal through the control horizon.
	\item apply the optimum control signal at time $k$ to the system, and get the updated state values.
	\item $k=k+1$, and go to the second step.
\end{itemize} 

Therefore, solving the optimization problem is complicated and time-consuming since the algorithm predicts the states and the inputs in each step, and calculates the cost function for the whole horizon. The more accurate the prediction model is, the less complicated optimization problem, and the faster controller will be attained. \\

\subsection{THE INDIRECT ADAPTIVE FUZZY MODEL PREDICTIVE CONTROL}
\begin{figure}[thpb]
	\begin{center}
		\framebox{\parbox{3in}{
				\centering
				\includegraphics[height=7cm, width=7cm]{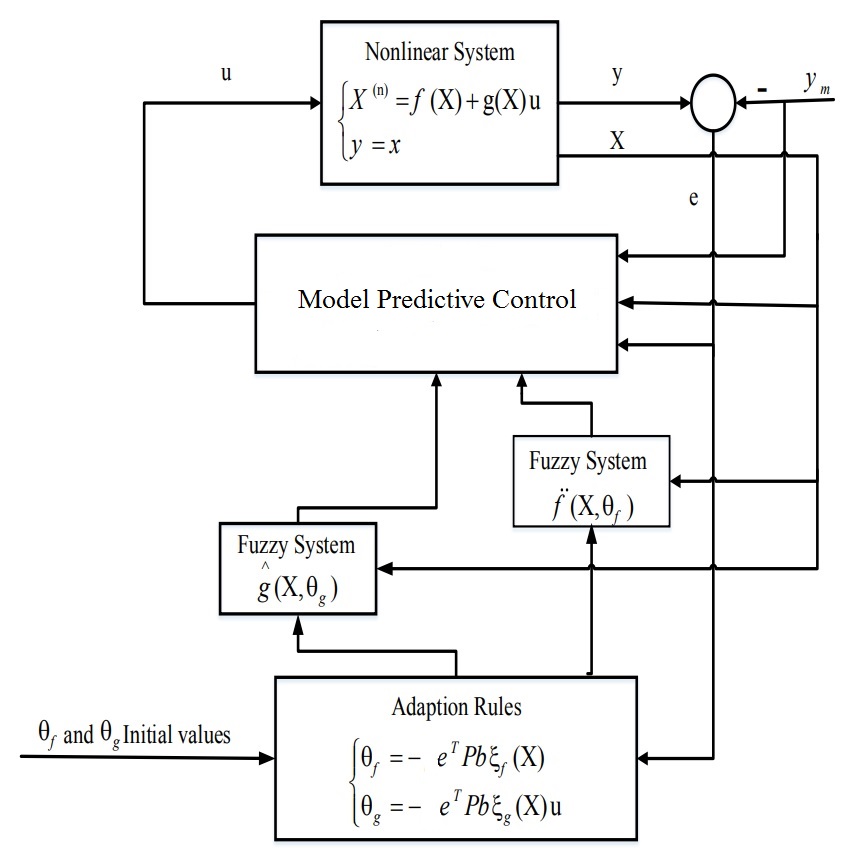}
		}}
		\caption{Block diagram of the indirect adaptive fuzzy model predictive control}
		\label{fuzzyMPC}
	\end{center}
\end{figure}

In the proposed model-based predictive control approach, an indirect adaptive fuzzy model is considered as the prediction model. The parameters of Mamdani fuzzy system are tuned based on the adaption rules. The rules are designed such that the closed-loop system stability is maintained (based on the Lyapunov criterion), and the system tracking error is minimized. \\ \indent
A nonlinear system can be linearized up to its relative degree $r$ as (\ref{linearized}) [8]. \\
\begin{fleqn}
	\begin{equation}  
	\begin{aligned}[b]
	&x^(r)=f(x,{x_1} \dot,\hdots,x^{(r-1)})+g(x,{x_1} \dot,\hdots,x^{(r-1)})u   
	\end{aligned}
	\label{linearized}
	\end{equation}
\end{fleqn}  \\ 
Above, $f$ and $g$ are the unknown functions. \\ \indent
Therefore, $\hat{f}$ and $\hat{g}$ are defined as the estimations of functions $f$ and $g$ respectively. If-then fuzzy rules are used to get the approximate functions. The fuzzy rules are provided from the input-output behavior of the system. Moreover, some parameters of $\hat{f}$ and $\hat{g}$ are free parameters and can vary during the process. Hence, the estimation functions are stated as (\ref{fhat&ghat}). \\
\begin{fleqn}
	\begin{equation}  
	\begin{aligned}[b]
	&\hat{f}=\hat{f}(X,\theta_f)\\
	&\hat{g}=\hat{g}(X,\theta_g)
	\end{aligned}
	\label{fhat&ghat}
	\end{equation}
\end{fleqn}  \\ 
where $\theta_f$ and $\theta_g$ are the free parameters of the estimation functions.  \\ \indent
Thus, based on the Mamdani fuzzy system and the system's input-output information, the estimation functions can be attained in a three-step process. In the first step, for each state variable $x_i$ of the system, $P_i$ number of fuzzy sets as ${A_i}^{l_i},\,\,\, l_i=1,2,\hdots,P_i$ are considered. Next, $q$ number of fuzzy sets as $B^l\,\,\, l=1,2,\hdots,q$ are considered for the system outputs. In the third step, the $\Pi_{i=1}^{n}P_i$ number of if-then fuzzy rules are constructed based on the system's input-output behavior as (\ref{ifthen_rule}) [8]. \\
\begin{fleqn}
	\begin{equation}  
	\begin{aligned}[b]
	&\text{if}\,\,\,x_i\,\,\,\text{is}\,\,\,{A_i}^{l_i},\,\,\,\text{then}\,\,\,\hat{f}(\hat{g})\,\,\,\text{is}\,\,\,B^{l_{i_1} \hdots l_{i_n}}\\
	&1\leq i_1,\hdots,i_n\leq n\,\,\,\text{and}\,\,\, 1\leq l_{i_1},\hdots,l_{i_n}\leq q
	\end{aligned}
	\label{ifthen_rule}
	\end{equation}
\end{fleqn}  \\ \indent
Introducing the $\Pi_{i=1}^{n}P_i$ dimension vectors $\varepsilon_f$ and $\varepsilon_g$ as (\ref{epsilon}); \\ 
\begin{fleqn}
	\begin{equation}  
	\begin{aligned}[b]
	&\varepsilon_f(X)= \frac{\Pi_{i=1}^{n}\mu _{A_i}^{l_i}(x_i)}{\sum_{l_1=1}^{P_1} \hdots \sum_{l_n=1}^{P_n}[\Pi_{i=1}^{n}\mu _{A_i}^{l_i}(x_i)]}\\
	&\varepsilon_g(X) = \frac{\Pi_{i=1}^{n}\mu _{A_i}^{l_i}(x_i)}{\sum_{l_1=1}^{P_1} \hdots \sum_{l_n=1}^{P_n}[\Pi_{i=1}^{n}\mu _{A_i}^{l_i}(x_i)]}
	\end{aligned}
	\label{epsilon}
	\end{equation}
\end{fleqn}  \\  
the estimation functions $\hat{f}$ and $\hat{g}$ are defined as (\ref{hatf&g}).  \\
\begin{fleqn}
	\begin{equation}  
	\begin{aligned}[b]
	&\hat{f}(X, \theta_f)= {\theta_f}^T\varepsilon_f(X)\\
	&\hat{g}(X, \theta_g)= {\theta_g}^T\varepsilon_g(X)
	\end{aligned}
	\label{hatf&g}
	\end{equation}
\end{fleqn}  \\  \indent
Above, g(X) and f(X) are actual functions from the state-space equation (\ref{eq:statespace}); $f(X)=a_2x_2+a_3sinx_3+a_4x_4$ and $g(X)=b_4$. Moreover, the membership functions $\mu_i$ are defined in the MATLAB fuzzy toolbox [8]. \\ \indent 
The adaption law is defined based on the Lyapunov theory. The Lyapunov function V is assumed as (\ref{lyapunov}). \\
\begin{fleqn}
	\begin{equation}  
	\begin{aligned}[b]
	&V = \frac{1}{2}e^TPe+e^TPbw + \frac{1}{2}(\theta _f-{\theta _f}^*)^T (\theta _f-{\theta _f}^*) \\
	&+ \frac{1}{2}(\theta _g-{\theta _g}^*)^T (\theta _g-{\theta _g}^*)  \\
	& A^TP+PA = -Q
	\end{aligned}
	\label{lyapunov}
	\end{equation}
\end{fleqn}  \\  \indent
Derivation of the Lyapunov function is defined as (\ref{lyapunov_derivation}). \\
\begin{fleqn}
	\begin{equation}  
	\begin{aligned}[b]
	&\dot V = -\frac{1}{2}e^TPe+e^TPbw + (\theta _f-{\theta _f}^*)^T [\dot {\theta_f} + e^TPb\varepsilon_f(X)] \\
	&+ (\theta _g-{\theta _g}^*)^T [\dot {\theta_g} + e^TPb\varepsilon_g(X)u]  
	\end{aligned}
	\label{lyapunov_derivation}
	\end{equation}
\end{fleqn}  \\  \indent
Matrix $Q$ is chosen as a positive definite matrix, and P is resulted by solving the Lyapunov equation in (\ref{lyapunov}). $e$ is the tracking error; the difference of the desired trajectory $y_m$ and the system's output $y$. Therefor, to attain the closed-loop stability, the Lyapunov derivative must be negative. So, the Lyapunov function's last two terms are forced to be zero. \\
\begin{fleqn}
	\begin{equation}  
	\begin{aligned}[b]
	&\dot {\theta _f} = - e^TPb\varepsilon_f(X) \\
	&\dot {\theta _g} =  - e^TPb\varepsilon_g(X)u
	\end{aligned}
	\label{lyapunov_stabled}
	\end{equation}
\end{fleqn}  \\  \indent
Thus, in each step of the proposed algorithm, the predicted signals from the fuzzy system are injected to the model predictive controller and the control law is generated. The current control signal is transferred to the system and the output signal is used for the next step. Therefore, the whole block diagram for the indirect adaptive fuzzy model predictive control strategy is shown in Fig. \ref{fuzzyMPC} [8-12].   \\
\\ 

\section{SIMULATION AND RESULTS}
\begin{figure}[thpb]
	\begin{center}
		\framebox{\parbox{3.5in}{
				\centering
				\includegraphics[height=5cm, width=9cm]{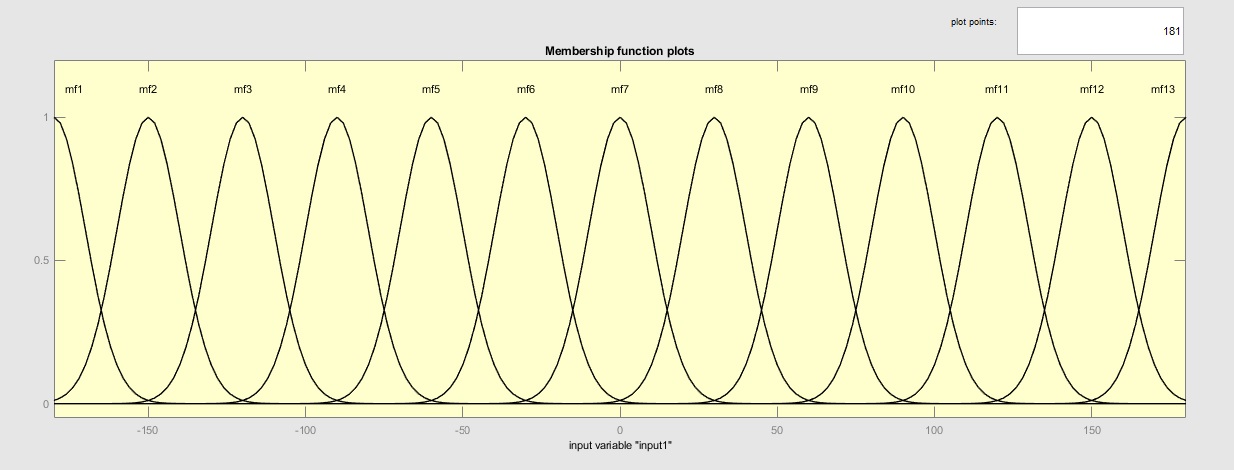}
		}}
		\caption{State variables' membership functions}
		\label{membership1}
	\end{center}
\end{figure}

\begin{figure}[thpb]
	\begin{center}
		\framebox{\parbox{3.5in}{
				\centering
				\includegraphics[height=5cm, width=9cm]{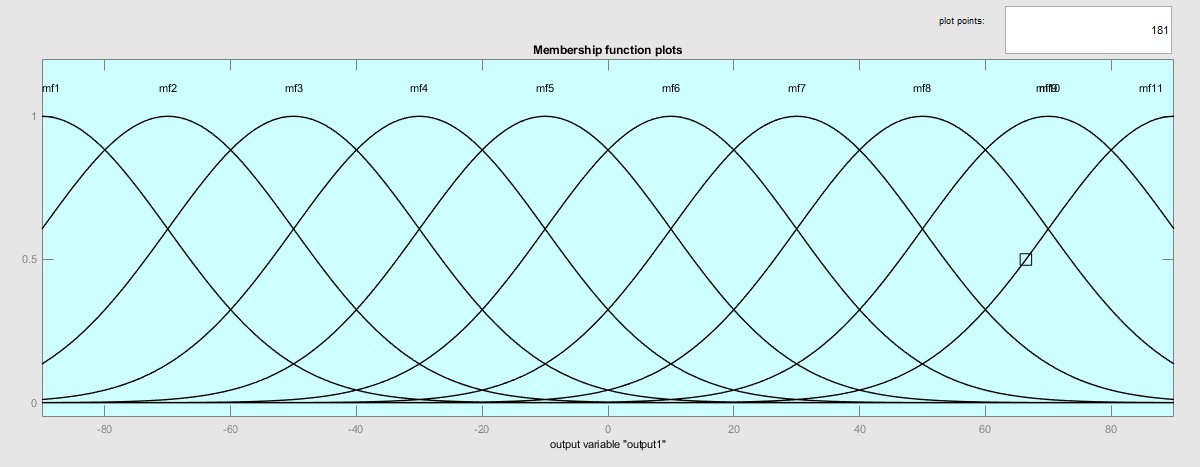}
		}}
		\caption{Block diagram of indirect adaptive fuzzy model predictive control}
		\label{membership2}
	\end{center}
\end{figure}

\begin{figure}[thpb]
	\begin{center}
		\framebox{\parbox{3.2in}{
				\centering
				\includegraphics[height=5cm, width=8cm]{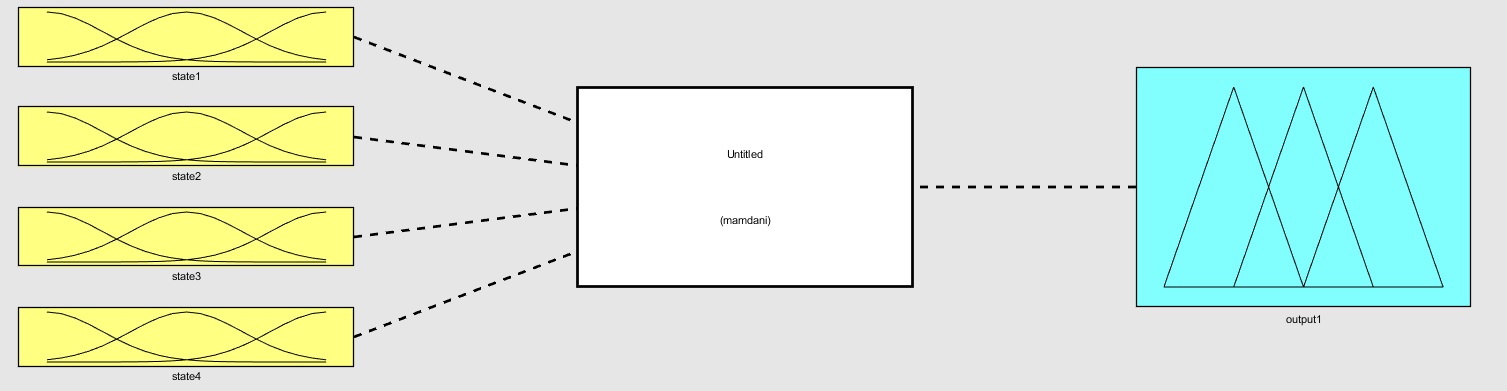}
		}}
		\caption{Fuzzy logic designer toolbox}
		\label{fuzz_toolbox}
	\end{center}
\end{figure}

\begin{figure}[thpb]
	\begin{center}
		\framebox{\parbox{3.3in}{
				\centering
				\includegraphics[height=6cm, width=8.5cm]{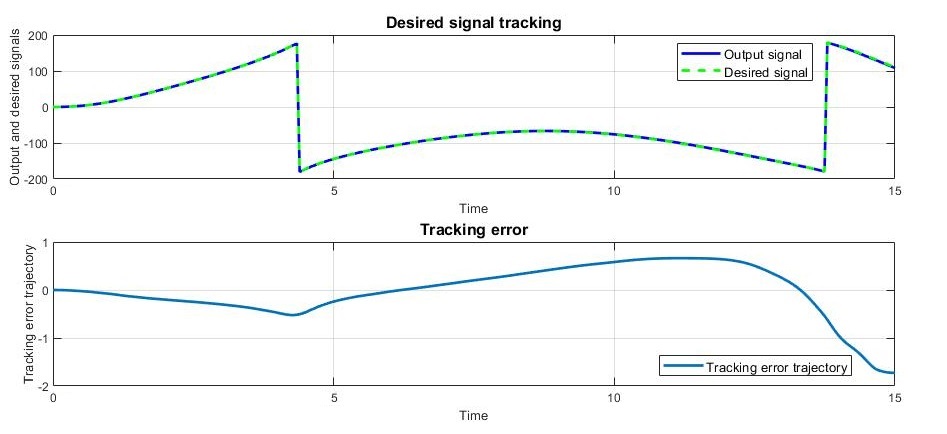}
		}}
		\caption{Desired signal tracking and tracking error using the classical MPC}
		\label{classic}
	\end{center}
\end{figure}

\begin{figure}[thpb]
	\begin{center}
		\framebox{\parbox{3.3in}{
				\centering
				\includegraphics[height=6cm, width=8.5cm]{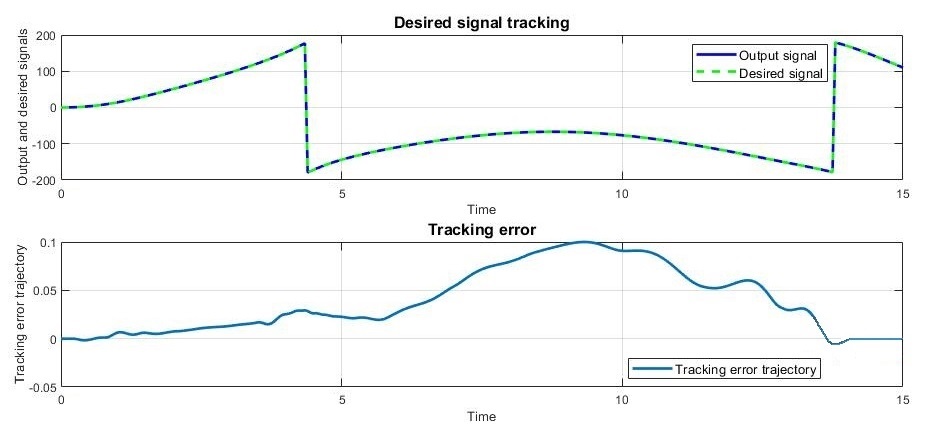}
		}}
		\caption{Desired signal tracking and tracking error using the indirect adaptive fuzzy MPC}
		\label{proposed}
	\end{center}
\end{figure}

Matrix A in the Lyapunov function (\ref{lyapunov}) is assumed as (\ref{matrixA}), and matrix Q of the Lyapunov equations is chosen with diagonal values $500$. \\
\begin{equation}      
{A}=
\begin{bmatrix} 
0 & 10 & 0 & 0 \\
0 & 0 & 10 & 0 \\
0 & 0 & 0 & 10 \\
-17.2 & -20.5 & -10 & 7 
\end{bmatrix}
\label{matrixA}
\end{equation}  \\  \indent
By solving the Lyapunov equation, matrix P is attained as (\ref{matrix_P}).  \\
\begin{fleqn}
	\begin{equation}  
	{P}=
	\begin{bmatrix} 
	75.5 & 35.1 & -12.2 & 18.3 \\
	33.7 & -5.8 & 9 & 10.1 \\
	-50 & -33 & 23.1 & -19.9 \\
	17.7 & -20.5 & -4.6 & -8.7 
	\end{bmatrix}
	\label{matrix_P}
	\end{equation}
\end{fleqn}  \\  \indent
The membership functions for the system's state variables and outputs are assumed as Figs. \ref{membership1} and \ref{membership2}. The membership functions are chosen as Gaussian functions. Besides, the membership functions are used in the if-then fuzzy rules based on Mamdani fuzzy system (Fig. \ref{fuzz_toolbox}).
  \\  \indent
In the model predictive control algorithm, the control horizon $K_c$ and the prediction horizon $K_p$ are chosen as 3 and 5 respectively. Moreover, the weight matrices Q and R are chosen as $0.1\times {I_{4\times4}}$ and $0.3$ respectively. Besides, the optimization problem is solved in each step through the fmincon optimization tool. \\  \indent
Figs. \ref{classic} and \ref{proposed} are resulted from the classical model predictive control approach and the proposed indirect adaptive fuzzy predictive control method respectively. \\  
  \\ \indent
According to the tracking error figures, the steady state error using the classical MPC strategy does not converge to zero, however, using the proposed indirect adaptive fuzzy predictive control, the error value throughout the process is significantly lower and eventually converges to zero. \\ \indent
Furthermore, the computation time for solving the optimization problem in the proposed strategy is lower than that of a classical MPC, because the initial prediction model is more accurate and there is no need to repeat the prediction process for a long time. The important feature of the proposed MPC algorithm is that it guarantees the closed-loop system stability through the Lyapunov theorem and it provides the steady-state zero tracking error for the controlled system. Besides, the proposed algorithm considers the uncertainties and disturbances prediction through the adaptive fuzzy approach which results in improved performance for the rotational inverted pendulum compared to the classical MPC. \\  \indent
To the best of our knowledge, the combination of the fuzzy modeling approach and model predictive control method is applied for the first time to the nonlinear rotational inverted pendulum benchmark problem.\\

\section{CONCLUSIONS}
A classical model predictive controller and an indirect adaptive fuzzy predictive controller were designed and implemented for a nonlinear rotational inverted pendulum with uncertainties. The nonlinear inverted pendulum system is a well-known benchmark problem in controls theory and applications. Choosing a prediction model for designing a predictive controller for the inverted pendulum is a prominent problem since the plant is highly nonlinear and subject to time-varying model uncertainties. \\ \indent
The proposed MPC strategy was used in combination with the indirect adaptive fuzzy model. A nonlinear fuzzy model was considered as the prediction model. The fuzzy system's parameters were then tuned based on the adaption rules such that the Lyapunov stability criterion was maintained and the system tracking error was minimized. The results of indirect adaptive fuzzy predictive control methodology were compared to the results of the classical model predictive control approach. Using the proposed method, the computation time for solving the optimization problem improved. Furthermore, the steady-state tracking error converged to zero using the proposed method, whereas the same value using the classical approach fluctuated between -2 and 1 throughout the whole simulation. It is worth mentioning that, combining the fuzzy modeling approach and model predictive control method was proposed for the first time for the nonlinear rotational inverted benchmark problem. Besides, the proposed strategy utilized the predictions of the model uncertainties and disturbances~ throughout the~MPC algorithm. 


\end{document}